\documentclass[number,citesort,seceqn,dvips]{arxbj}
\usepackage{upgreek}

\aid{0}
\volume{17}
\issue{2}
\pubyear{2011}
\firstpage{507}
\lastpage{529}
\doi{10.3150/10-BEJ282}

\makeatletter
\newcommand{\E}{\mathbb{E}}
\newcommand{\Z}{\mathbb{Z}}

\newcommand{\R}{\mathbb{R}}
\newcommand{\sR}{\mathbb{R}}
\newcommand{\Var}{\operatorname{Var}}
\newcommand{\lead}{\operatorname{lead}}
\newcommand{\Cov}{\operatorname{Cov}}

\newcommand{\D}{\mathbb{D}}

\newtheorem{theo}{Theorem}[section]
\newtheorem{lem}{Lemma}[section]

\newremark{example}{Example}[section]
\newremark{application}{Application}[section]
\makeatother

\begin{document}
\begin{frontmatter}

\title{An extended Stein-type covariance identity for the Pearson family with applications
to lower variance bounds}
\runtitle{An extended Stein-type identity}

\begin{aug}
\author[1]{\fnms{G.} \snm{Afendras}\ead[label=e1]{afendras.giorgos@ucy.ac.cy}\thanksref{1}},
\author[2]{\fnms{N.} \snm{Papadatos}\corref{}\ead[label=e2,mark]{npapadat@math.uoa.gr}\thanksref{2,e2}}
\and
\author[2]{\fnms{V.} \snm{Papathanasiou}\ead[label=e3,mark]{bpapath@math.uoa.gr}\thanksref{2,e3}}

\runauthor{G. Afendras, N. Papadatos and V. Papathanasiou}

\address[1]{Department of Mathematics and Statistics,
University of Cyprus,
P.O. Box 20537,
1678 Nicosia, Cyprus.
\printead{e1}}

\address[2]{Department of Mathematics,
Section of Statistics and O.R.,
University of Athens,
Panepistemiopolis,
157 84 Athens, Greece.
\printead{e2,e3}}
\end{aug}

\received{\smonth{6} \syear{2008}}
\revised{\smonth{2} \syear{2010}}

%
\begin{abstract}
For an absolutely continuous (integer-valued)
r.v. $X$ of the Pearson (Ord) family,
we show
that, under natural moment conditions, a Stein-type covariance
identity of order $k$ holds (cf. [Goldstein and Reinert, \textit{J. Theoret. Probab.}~\textbf
{18} (2005) 237--260]).
This identity is closely related to the
corresponding sequence of orthogonal polynomials, obtained by a
Rodrigues-type formula, and provides convenient
expressions for the Fourier coefficients of an arbitrary function.
Application of the covariance identity yields some novel
expressions for the corresponding lower variance bounds for a
function of the r.v. $X$, expressions that seem to be known only in
particular cases (for the Normal, see [Houdr\'{e} and Kagan, \textit{J. Theoret. Probab.}
\textbf{8} (1995)
23--30]; see also
[Houdr\'{e} and P\'{e}rez-Abreu, \textit{Ann. Probab.} \textbf{23}
(1995) 400--419]
for corresponding results related to the Wiener and
Poisson processes).
Some applications are also given.
\end{abstract}

%

\end{frontmatter}
%
\section{Introduction}
\label{sec1}

For an r.v. $X$ with density $f$, mean $\mu$ and finite variance
$\sigma^2$, Goldstein and Reinert \cite{18}
showed
the identity (see also \cite{26})
\begin{equation}\label{eq1.1}
\Cov(X,g(X))=\sigma^2\E[g'(X^*)],
\end{equation}
which holds for any absolutely continuous function
$g\dvtx \R\to\R$ with a.s. derivative $g'$ such that the right-hand side
is finite.
In (\ref{eq1.1}), $X^*$ is defined to be the r.v. with density
$f^*(x)=\frac{1}{\sigma^2}
\int_{-\infty}^x
(\mu-t)f(t)\,\mathrm{d}t=\frac{1}{\sigma^2}\int_x^{\infty}(t-\mu)f(t)\,\mathrm{d}t$,
$x\in\R$.

Identity (\ref{eq1.1}) extends the well-known Stein
identity
for the standard normal \cite{33,34};
a discrete
version of (\ref{eq1.1}) can be found in, for example,
\cite{13},
where the derivative has been replaced by
the forward difference of $g$. In particular, identities
of the form (\ref{eq1.1}) have
many applications
to variance bounds and
characterizations \cite{4,13,26},
and to approximation
procedures \cite{11,14,15,18,27,31,33}.
Several extensions and applications can be found in \cite{10,19,29}.

In \cite{23}, the (continuous) Pearson family is parametrized
by the fact that
there exists a quadratic $q(x)=\delta
x^2+\beta x+\gamma$ such that
\begin{equation}\label{eq1.2}
\int_{-\infty}^x (\mu-t)f(t) \,\mathrm{d}t = q(x)f(x),\qquad x\in\R.
\end{equation}
Typically, the usual definition of a Pearson r.v. is related to the
differential equation
$f'(x)/f(x)=(\alpha-x)/p_2(x)$, with $p_2$ being a polynomial of
degree at most 2.
In fact, the set-up of (\ref{eq1.2}) (including, e.g., the standard uniform
distribution with $q(x)=x(1-x)/2$) will be the framework of the
present work and will hereafter be called ``the Pearson
family of continuous distributions''. It is easily seen that under
(\ref{eq1.2}), the support of $X$, $S(X)=\{x\dvtx f(x)>0\}$, must be an interval,
say $(r,s)$ with $-\infty\leq r <s\le\infty$, and $q(x)$ remains
strictly positive for $x\in(r,s)$. Clearly, under (\ref{eq1.2}), the
covariance identity (\ref{eq1.1}) can be rewritten as
\begin{equation}\label{eq1.3}
\E
[(X-\mu)g(X)]=\E[q(X)g'(X)].
\end{equation}
It is known that,
under appropriate
moment conditions, the functions
\begin{equation}\label{eq1.4}
P_k(x)=\frac{(-1)^k}{f(x)}\frac{\mathrm{d}^k}{\mathrm{d} x^k} [q^k(x)f(x)],\qquad x\in(r,s),\ k=0,1,\ldots,M
\end{equation}
(where $M$ can be finite or infinite) are orthogonal polynomials
with respect to the density $f$ so that the quadratic $q(x)$ in
(\ref{eq1.2}) generates a sequence of orthogonal polynomials by the
Rodrigues-type formula (\ref{eq1.4}).
In fact, this approach is related to the Sturm--Liouville
theory, \cite{17}, Section~5.2; see also \cite{24,28}.

In the present paper, we provide an extended Stein-type identity of
order $k$
for the Pearson family. This identity takes
the form
\begin{equation}\label{eq1.5}
\E[ P_k(X)g(X)]= \E\bigl[ q^k(X)g^{(k)}(X)\bigr],
\end{equation}
where
$g^{(k)}$ is the $k$th derivative of $g$ (since $P_1(x)=x-\mu$,
(\ref{eq1.5}) for $k=1$ reduces to (\ref{eq1.3})). Identity (\ref{eq1.5}) provides a
convenient
formula for the $k$th Fourier coefficient of $g$, corresponding to
the orthogonal polynomial $P_k$ in (\ref{eq1.4}). For its proof, we
make use of a novel ``Rodrigues inversion'' formula that may be of
some
interest in itself.
An identity similar to (\ref{eq1.5}) holds for the discrete Pearson (Ord)
family.
Application of (\ref{eq1.5}) and its discrete analog yields the
corresponding lower variance bounds, obtained in Section \ref{sec4}. The
lower bound for the $\operatorname{Poisson}(\lambda)$ distribution, namely,
\begin{equation}\label{eq1.6}
\Var g(X)\geq\sum_{k=1}^n \frac{\lambda^{k}}{k!}\E^2
[\Delta^{k}[g(X)]]
\end{equation}
(cf. \cite{21})
and the
corresponding one for the Normal$(\mu,\sigma^2)$ distribution
\cite{20},
\begin{equation}\label{eq1.7}
\Var g(X)\geq\sum_{k=1}^n \frac{(\sigma^2)^{k}}{k!}\E^2
\bigl[g^{(k)}(X)\bigr],
\end{equation}
are particular examples (Examples \ref{ex4.1} and \ref{ex4.5}) of Theorems \ref{th4.1} and
\ref{th4.2}, respectively.
Both (\ref{eq1.6}) and~(\ref{eq1.7}) are particular cases of the finite form
of Bessel's inequality and, under completeness, they can
be extended to the corresponding Parseval identity.
In Section \ref{sec5}, we
show that this can be done
for a fairly large family of r.v.'s, including, of course,
the normal, the Poisson and, in general,
all the r.v.'s of the Pearson system
which have finite moments of any order.
For instance, when $X$ is $\operatorname{Normal}(\mu,\sigma^2)$,
inequality (\ref{eq1.7}) (and identity (\ref{eq1.5}))
can be strengthened to the covariance identity
\begin{equation}\label{eq1.8}
\Cov[g_1(X),g_2(X)]= \sum_{k=1}^{\infty} \frac{(\sigma^2)^{k}}{k!}
\E\bigl[g_1^{(k)}(X)\bigr] \E\bigl[g_2^{(k)}(X)\bigr],
\end{equation}
provided that for $i=1,2$, $g_i\in\D^{\infty}(\R)$,
$\E|g_i^{(k)}(X)|<\infty$, $k=1,2,\ldots,$
and that $\E[g_i(X)]^2<\infty$. Similar identities hold
for Poisson, negative binomial, beta and gamma
distributions.
These kinds of variance/covariance expressions
may sometimes be useful in inference problems -- see, e.g., the Applications \ref{ap5.1} and \ref{ap5.2} at
the end of the paper.

\section{Discrete orthogonal polynomials and the covariance
identity}\label{sec2}

In order to
simplify notation, we assume that $X$ is a non-negative
integer-valued r.v. with mean $\mu<\infty$.
We also assume that
there exists a quadratic $q(x)=\delta x^2+\beta x+\gamma$
such that
\begin{equation}\label{eq2.1}
\sum_{j=0}^x(\mu-j)p(j)=q(x)p(x),\qquad x=0,1,\ldots,
\end{equation}
where $p(x)$ is the probability
function of $X$. Relation (\ref{eq2.1}) describes the
discrete Pearson system (Ord family) \cite{23}. Let $\Delta^k$
be the forward difference operator defined by
$\Delta[g(x)]=g(x+1)-g(x)$ and $\Delta^k[g(x)]=\Delta[\Delta^{k-1}[g(x)]]$
($\Delta^0[g]\equiv g$, $\Delta^1\equiv\Delta$).
We also set $q^{[k]}(x)=q(x)q(x+1)\cdots q(x+k-1)$ (with
$q^{[0]}\equiv1$, $q^{[1]}\equiv q$).

We first show some useful lemmas.

\begin{lem}\label{le2.1}
If $h(x)=0$ for $x<0$ and
%
%
\begin{eqnarray}\label{eq2.2}
\sum_{x=0}^{\infty} |\Delta^j [h(x-j)] \Delta^{k-j}[g(x)]|&<&\infty
\qquad\mbox{for }
j=0,1,\ldots,k,
\\
\label{eq2.3}\lim_{x\to\infty} \Delta^j [h(x-j)] \Delta^{k-j-1}[g(x)]&=&0\qquad
\mbox{for } j=0,1,\ldots,k-1,
\end{eqnarray}
then
\begin{eqnarray}\label{eq2.4}
&&(-1)^k \sum_{x=0}^{\infty} \Delta^k[h(x-k)]g(x)
\nonumber
\\[-8pt]
\\[-8pt]
\nonumber
&&\quad=
\sum_{x=0}^{\infty} h(x)\Delta^k[g(x)].
\end{eqnarray}
\end{lem}
\begin{pf}
We have
\begin{eqnarray*}
\sum_{x=0}^{\infty} h(x) \Delta^k[g(x)]
&=&
\lim_{n\to\infty} \sum_{x=0}^n h(x) \bigl(\Delta^{k-1}[g(x+1)]-\Delta
^{k-1}[g(x)]\bigr)
\\
&=&
\lim_{n\to\infty} \Biggl[ h(n+1)\Delta^{k-1}[g(n+1)]-\sum_{x=0}^{n+1}
\Delta[h(x-1)] \Delta^{k-1}[g(x)] \Biggr]
\\
&=&
\lim_{n\to\infty}
h(n+1)\Delta^{k-1}[g(n+1)]-\sum_{x=0}^{\infty} \Delta[h(x-1)]
\Delta^{k-1}[g(x)]
\\
&=&
-\sum_{x=0}^{\infty} \Delta[h(x-1)]
\Delta^{k-1}[g(x)].
\end{eqnarray*}
By the same calculation, it follows that
\[
(-1)^j \sum_{x=0}^{\infty} \Delta^j [h(x-j)]
\Delta^{k-j}[g(x)]=(-1)^{j+1} \sum_{x=0}^{\infty} \Delta^{j+1}
[h(x-j-1)] \Delta^{k-j-1}[g(x)]
\]
for any $j\in\{0,1,\ldots,k-1\}$.
\end{pf}

\begin{lem}\label{le2.2} For each $n\geq0$, there exist polynomials $Q_{i,n}(x)$,
$i=0,1,\ldots,n$, such that the degree of each $Q_{i,n}$ is at
most $i$ and
\begin{eqnarray}\label{eq2.5}
\Delta^i \bigl[q^{[n]}(x-n)p(x-n)\bigr]
=q^{[n-i]}(x-n+i)p(x-n+i)Q_{i,n}(x),\qquad
i=0,1,\ldots, n.\hspace*{26pt}
\end{eqnarray}
Moreover, the leading coefficient (i.e., the coefficient of
$x^n$) of $Q_{n,n}$ is given by $\lead(Q_{n,n})=(-1)^n
\prod_{j=n-1}^{2n-2}(1-j\delta)$, where an empty product should be
treated as $1$.
\end{lem}
\begin{pf}
For $n=0$, the assertion is obvious and $Q_{0,0}(x)=1$. For $n=1$,
the assertion follows from the assumption (\ref{eq2.1}) with
$Q_{0,1}(x)=1$, $Q_{1,1}(x)=\mu-x$. For the case $n\geq2$, the
assertion will be proven using (finite) induction on $i$. Indeed,
for $i=0$, (\ref{eq2.5}) holds with $Q_{0,n}(x)=1$. Assuming that the
assertion holds for some $i\in\{0,1,\ldots, n-1\}$ and setting
$h_n(x)=q^{[n]}(x-n)p(x-n)$, it follows that
\begin{eqnarray*}
\Delta^{i+1} h_n(x)
&=&
\Delta[ \Delta^i h_n (x)]
\\
&=&
\Delta[h_{n-i}(x) Q_{i,n}(x)]
\\
&=&
\Delta\bigl[
q(x-n+i)p(x-n+i)\bigl(q^{[n-i-1]}(x-n+i+1)Q_{i,n}(x)\bigr)\bigr]
\\
&=&
q(x-n+i+1)p(x-n+i+1) \Delta\bigl[q^{[n-i-1]}(x-n+i+1)Q_{i,n}(x)\bigr]
\\
&&{}
+\bigl(\mu-(x-n+i+1)\bigr)p(x-n+i+1)q^{[n-i-1]}(x-n+i+1)Q_{i,n}(x),
\end{eqnarray*}
where the obvious relation $\Delta[q(x)p(x)]=(\mu-(x+1))p(x+1)$,
equivalent to (\ref{eq2.1}), has been used with $x-n+i$ in place of $x$.
Moreover,
\begin{eqnarray*}
&&\Delta\bigl[q^{[n-i-1]}(x-n+i+1)Q_{i,n}(x)\bigr]\\
&&\quad=
q^{[n-i-1]}(x-n+i+2)\Delta[Q_{i,n}(x)]
+
Q_{i,n}(x)\Delta\bigl[q^{[n-i-1]}(x-n+i+1)\bigr]
\end{eqnarray*}
and
\[
\Delta\bigl[q^{[n-i-1]}(x-n+i+1)\bigr]
=
q^{[n-i-2]}(x-n+i+2)
\bigl(q(x)-q(x-n+i+1)\bigr)
\]
(where, for $i=n-1$, the right-hand side of the above should be treated as
$0$). Observing that
\begin{eqnarray*}
&&q(x-n+i+1)q^{[n-i-2]}(x-n+i+2)\bigl(q(x)-q(x-n+i+1)\bigr)
\\
&&\quad =q^{[n-i-1]}(x-n+i+1)\bigl(q(x)-q(x-n+i+1)\bigr)
\end{eqnarray*}
(which is also true in the case where $i=n-1$), the above
calculations show that (\ref{eq2.5}) holds with
\[
Q_{i+1,n}(x)= P_{i,n}(x) Q_{i,n}(x)+ q(x)\Delta[Q_{i,n}(x)],
\]
where $P_{i,n}(x)=\mu-(x-n+i+1)+q(x)-q(x-n+i+1)$ is a linear
polynomial or a constant. From the above recurrence, it follows
immediately that
$\lead(Q_{i+1:n})=-(1-(2n-i-2)\delta)\lead(Q_{i,n})$,
$i=0,1,\ldots,n-1$; this, combined with the fact that
$\lead(Q_{0,n})=1$, yields the desired result.
\end{pf}

It is easy to see that under (\ref{eq2.1}), the support of $X$,
$S(X)=\{x\in\Z\dvtx  p(x)>0\}$, is a finite or infinite integer
interval. This integer interval will be denoted by $J$.
Here, the term ``integer interval'' means that ``if $j_1$ and~$j_2$ are integers belonging to $J$, then all integers
between $j_1$ and~$j_2$ also belong to $J$''.

\begin{lem}\label{le2.3}
For each $k=0,1,2,\ldots,$ define the functions $P_k(x)$, $x\in J$,
by the Rodrigues-type formula
\begin{eqnarray}\label{eq2.6}
P_k(x)&=&\frac{(-1)^k}{p(x)} \Delta^k
\bigl[q^{[k]}(x-k)p(x-k)\bigr]
\nonumber
\\[-8pt]
\\[-8pt]
\nonumber
&=&\frac{1}{p(x)}\sum_{j=0}^k
(-1)^{k-j}{k\choose j}q^{[k]}(x-j)p(x-j).
\end{eqnarray}
We then have the following:
\begin{longlist}[(a)]
\item[(a)] Each $P_k$ is a polynomial of degree at most $k$, with
$\lead(P_k)=\prod_{j=k-1}^{2k-2}(1-j\delta)$ \textup{(in the sense
that the function $P_k(x)$, $x\in J$, is the restriction of a real
polynomial $G_k(x)=\sum_{j=0}^k c(k,j)x^j$, $x\in\R$, of degree at
most $k$, such that $c(k,k)=\lead(P_k)$)}.\vspace*{2pt}

\item[(b)] Let $g$ be an arbitrary function defined on $J$
\textup{(the integer interval support of $X$)} and, if $J$ is infinite,
assume in addition that the functions $g$ and $h(x):=q^{[k]}(x)p(x)$
satisfy the requirements~\textup{(\ref{eq2.2})} and \textup{(\ref{eq2.3})} of Lemma
\ref{le2.1}.
Then,
\[
\E|P_{k}(X)g(X)|<\infty, \qquad \E \bigl[q^{[k]}(X)|\Delta^{k}[g(X)]|\bigr]<\infty
\]
and the following identity
holds:
\begin{equation}\label{eq2.7}
\E[P_k(X)g(X)]=\E\bigl[q^{[k]}(X)\Delta^{k}[g(X)]\bigr].
\end{equation}
\end{longlist}
\end{lem}
\begin{pf}
Part (a) follows from Lemma \ref{le2.2} since $P_k=(-1)^k Q_{k,k}$. For
part (b), assume first that $p(0)>0$, that is, the support $J$ is
either $J=\{0,1,\ldots\}$ or $J=\{0,1,\ldots,N\}$ for a positive
natural number $N$. In the unbounded case, (\ref{eq2.7}) follows by an
application of Lemma \ref{le2.1} to the functions $h(x)=q^{[k]}(x)p(x)$
and $g$ since they satisfy the conditions (\ref{eq2.2}) and (\ref{eq2.3}). For
the bounded case, it follows by (\ref{eq2.1}) that $x=N$ is a zero of
$q(x)$ so that $q(x)=(N-x)(\mu/N-\delta x)$, where, necessarily,
$\delta<\mu/(N(N-1))$. Thus, in the case where $k> N$,
$q^{[k]}(x)=0$ for all $x\in J$ so that the right-hand side of (\ref{eq2.7})
vanishes. On the other hand, the left-hand side of (\ref{eq2.7}) (with $P_k$
given by
(\ref{eq2.6})) is $\sum_{x=0}^N g(x) \sum_{j=0}^k (-1)^{k-j}{k\choose j}
q^{[k]}(x-j)p(x-j)$ and it is easy to verify that for all $x\in J
$ and all $j\in\{0,1,\ldots,k\}$, the quantity
$q^{[k]}(x-j)p(x-j)$ vanishes. Thus, when $k>N$, (\ref{eq2.7}) holds in
the trivial sense $0=0$. For $k\leq N$,
the left-hand side of (\ref{eq2.7}) equals
\begin{eqnarray*}
&&\sum_{j=0}^k(-1)^{k-j} \pmatrix{k\cr j} \sum_{x=-j}^{N-j} q^{[k]}(x)
p(x)g(x+j)\\
&&\quad=
\sum_{j=0}^k(-1)^{k-j} \pmatrix{k\cr j}
\sum_{x=0}^{N-j} q^{[k]}(x)p(x)g(x+j)
\\
&&\quad=
\sum_{x=0}^{N} q^{[k]}(x)
p(x)\sum_{j=0}^{\min\{k,N-x\}}(-1)^{k-j} \pmatrix{k\cr
j}g(x+j)
\\
&&\quad=
\sum_{x=0}^{N-k} q^{[k]}(x) p(x) \sum_{j=0}^k (-1)^{k-j}
\pmatrix{k\cr j}g(x+j),
\end{eqnarray*}
where we have made use of the facts that $p(x)=0$ for $x<0$ and $q^{[k]}(x)=0$
for $N-k<x\leq N$. Thus, (\ref{eq2.7}) is proved under the assumption
$p(0)>0$. For the
general case, where the support of~$X$ is either $\{m,m+1,
\ldots,N\}$ or $\{m,m+1, \ldots\}$, it suffices to apply the same
arguments to the r.v. $X-m$.
\end{pf}

\begin{theo}\label{th2.1}
Suppose that $X$ satisfies \textup{(\ref{eq2.1})} and has $2n$ finite
moments
for some $n\geq1$.
The polynomials $P_k$, $k=0,1,\ldots,n$, defined by \textup{(\ref{eq2.6})},
then satisfy the orthogonality condition
\begin{equation}\label{eq2.8}
\E[P_k(X)P_m(X)]=\delta_{k,m}
k!\E\bigl[q^{[k]}(X)\bigr]\prod_{j=k-1}^{2k-2}(1-j\delta),\qquad
k,m=0,1,\ldots,n,
\end{equation}
where $\delta_{k,m}$ is Kronecker's delta.
\end{theo}
\begin{pf}
Let $0\leq m\leq k\leq n$. First, observe that $q^{[k]}$ is a
polynomial of degree at most $2k$ and~$P_m$ a polynomial of
degree at most $m$. Note that the desired result will be deduced
if we can apply~(\ref{eq2.7}) to the function $g=P_m$; this can be
trivially applied when $J$ is finite. In the remaining case where
$J$ is infinite, we have to verify conditions~(\ref{eq2.2}) and (\ref{eq2.3}) of
Lemma~\ref{le2.1} for the functions $h=q^{[k]}p$ and $g=P_m$,
or, more generally, for any polynomial $g$ of degree less than or
equal to $k$. Since the case $m=0$ is obvious ($P_0(x)=1$), we
assume that $1\leq m\leq k\leq n$. Observe that~(\ref{eq2.2}) is satisfied
in this case because of the assumption
$\E|X|^{2k}<\infty$.
Indeed, from Lemma~\ref{le2.2},
$\Delta^j[q^{[k]}(x-j)p(x-j)]=p(x)
q^{[k-j]}(x)Q_{j,k}(x-j+k)$,
where $Q_{j,k}$ is of degree at most
$j$ and $q^{[k-j]}$ is of degree at most $2k-2j$ so that their
product is a polynomial of degree at most $2k-j$, while
$\Delta^{k-j}[g]$ is a polynomial of degree at most $j$. On the other
hand, since $q(x)p(x)$ is eventually decreasing and
$\sum_{x=0}^{\infty} x^{2k-2} q(x)p(x)<\infty$, we have, as
$y\to\infty$, that
\begin{eqnarray*}
q(2y)p(2y)y(y+1)^{2k-2}&\leq&\sum_{x=y+1}^{2y} x^{2k-2}q(x)p(x)\to
0 \quad \mbox{and}
\\
q(2y+1)p(2y+1)(y+1)^{2k-1}&\leq&\sum_{x=y+1}^{2y+1}
x^{2k-2}q(x)p(x)\to0.
\end{eqnarray*}
Hence,
$\lim_{x\to\infty} x^{2k-1}q(x)p(x)=0$.
For $j\leq k-1$, we have, from Lemma \ref{le2.2}, that
\[
\Delta^j[h(x-j)]
\Delta^{k-j-1}[g(x)]=q(x)p(x)R(x),
\]
where
\[
R(x)=q^{[k-j-1]}(x+1)Q_{j,k}(x-j+k)\Delta^{k-j-1}[g(x)]
\]
is a
polynomial of degree at most $2k-1$. This, combined with the above
limit, verifies (\ref{eq2.3}) and a final application of (\ref{eq2.7}) completes the
proof.
\end{pf}

The covariance identity (\ref{eq2.7}) enables the calculation of the Fourier
coefficients of any function~$g$ in terms of its differences
$\Delta^{k}[g]$, provided that conditions (\ref{eq2.2}) and (\ref{eq2.3}) are fulfilled
for $g$ and $h(x)=q^{[k]}(x)p(x)$. Since this identity is important
for the applications, we state and prove a more general result that
relaxes these conditions. The proof of this result contains a novel
inversion formula for the polynomials obtained by the Rodrigues-type
formula (\ref{eq2.6}).

\begin{theo}\label{th2.2} Assume that $X$ satisfies \textup{(\ref{eq2.1})} and has $2k$ finite moments,
and suppose that for some function $g$ defined on $J$,
\[
\E
\bigl[
q^{[k]}(X) |\Delta^{k}[g(X)]|\bigr]<\infty.
\]
Then,
\[
\E|P_k(X)g(X)|<\infty
\]
and the
covariance identity \textup{(\ref{eq2.7})} holds. Moreover, provided that
$k\geq1$ and $X$ has $2k-1$ finite moments, the following
``Rodrigues inversion'' formula holds:
\begin{equation}\label{eq2.9}
q^{[k]}(x)p(x)
=\frac{1}{(k-1)!}\sum_{y=x+1}^{\infty}
(y-x-1)_{k-1}P_k(y) p(y),\qquad x=0, 1,\ldots,
\end{equation}
where, for
$x\in\R$, $(x)_n=x(x-1)\cdots(x-n+1)$, $n=1,2,\ldots,$ and $(x)_0=1$.
\end{theo}
\begin{pf}Relation (\ref{eq2.7}) is obvious when $k=0$. Also, the
assertion follows trivially by Lemma~\ref{le2.3}(b) if $J$ is finite.
Assume, now, that $k\geq1$ and that $J$ is unbounded, fix
$s\in\{0,1,\ldots\}$ and consider the function
\[
g_{s,k}(x):=\frac{1}{(k-1)!}I(x\geq s+1) (x-s-1)_{k-1},\qquad
x=0,1,\ldots.
\]
It is easily seen that the pair of functions $g=g_{s,k}$ and
$h(x)=q^{[k]}(x)p(x)$ satisfy (\ref{eq2.2}) and~(\ref{eq2.3}) and also that
$\Delta^{k}[g_{s,k}(x)]=I(x=s)$. Thus, by (\ref{eq2.7}), we get
(\ref{eq2.9}) (cf. formulae (\ref{eq3.6}) and (\ref{eq3.7}), below, for the continuous
case),
where $P_k$ is defined by (\ref{eq2.6}) and the series converges from
Lemma~\ref{le2.3}(a) and the fact that $X$ has $2k-1$ finite moments (note
that $2k-1$ finite moments suffice for this inversion formula).
Since $P_k$ is a polynomial and the left-hand side of (\ref{eq2.9}) is strictly
positive for all large enough $x$, it follows that $P_k(y)>0$ for
all large enough $y$. Using~(\ref{eq2.9}), the assumption on $\Delta^{k}[g]$ is
equivalent to the fact that
\[
\frac{1}{(k-1)!}\sum_{x=0}^\infty|\Delta^{k}[g(x)]|
\sum_{y=x+1}^{\infty} (y-x-1)_{k-1}P_k(y) p(y)<\infty
\]
and arguments similar to those used in the proof of Theorem
\ref{th3.1}(b) below show that
\[
\frac{1}{(k-1)!}\sum_{x=0}^\infty|\Delta^{k}[g(x)]|
\sum_{y=x+1}^{\infty} (y-x-1)_{k-1}|P_k(y)| p(y)<\infty.
\]
Therefore, we can interchange the order of summation, obtaining
\begin{eqnarray*}
\E\bigl[ q^{[k]}(X)\Delta^{k}[g(X)]\bigr]
&=&
\frac{1}{(k-1)!}
\sum_{x=0}^\infty
\Delta^{k}[g(x)] \sum_{y=x+1}^{\infty} (y-x-1)_{k-1}P_k(y) p(y)
\\
&=&
\frac{1}{(k-1)!}\sum_{y=1}^\infty P_k(y) p(y)
\sum_{x=0}^{y-1} (y-x-1)_{k-1}
\Delta^{k}[g(x)]
\\
&=&
\E[P_k(X) G(X)],
\end{eqnarray*}
where
\[
G(x)=\frac{1}{(k-1)!} \sum_{y=0}^{x-k} (x-1-y)_{k-1}\Delta
^{k}[g(y)],\qquad  x=0,1,\ldots,
\]
and where an empty sum should be treated as $0$. Taking forward
differences, it follows that $\Delta^{k}[G(x)]=\Delta^{k}[g(x)]$ so that
$G=g+H_{k-1}$, where $H_{k-1}$ is a polynomial of degree at most
$k-1$, and the desired result follows from the orthogonality of
$P_k$ to polynomials of degree lower than $k$.
This completes the proof.
\end{pf}

\section{The generalized Stein-type identity for the continuous
case}\label{sec3}

The orthogonality of polynomials (\ref{eq1.4})
has been shown, for example, in \cite{17};
see also \cite{24,28}.
For our purposes, we review some details.

The induction formula, as in (\ref{eq2.5}), here takes the form
\begin{equation}\label{eq3.1}
\frac{\mathrm{d}^{i}}{\mathrm{d}x^{i}}[q^k(x)f(x)]=q^{k-i}(x)f(x)Q_{i,k}(x),\qquad i=0,1,\ldots,k,
\end{equation}
with $Q_{0,k}(x)=1$, where $Q_{i,k}$ is a polynomial of degree at
most $i$ and the explicit recurrence for $Q_{i,k}$ is
\[
Q_{i+1,k}(x)=\bigl(\mu-x+(k-i-1)(2\delta x+\beta)\bigr)Q_{i,k}(x)+q(x)
Q'_{i,k}(x),\qquad i=0,1,\ldots,k-1.
\]
This immediately implies that
$\lead(P_k)=(-1)^k\lead(Q_{k,k})=\prod_{j=k-1}^{2k-2}(1-j\delta)$,
as in the discrete case; see also \cite{7,17}.
The covariance identity, as in (\ref{eq2.7}),
here takes the form (after repeated integration by parts; cf. \cite{22,28})
\begin{equation}\label{eq3.2}
\E[P_k(X)g(X)] =\E\bigl[q^k(X)g^{(k)}(X)\bigr],
\end{equation}
provided that the expectations are finite and
%
%
%
\begin{eqnarray}\label{eq3.3}
&&\lim_{x\to r+}q^{i+1}(x)f(x)Q_{k-i-1,k}(x) g^{(i)}(x)
\nonumber
\\[-8pt]
\\[-8pt]
\nonumber
&&\quad=\lim_{x\to
s-}q^{i+1}(x)f(x)Q_{k-i-1,k}(x) g^{(i)}(x)=0,\qquad
i=0,1,\ldots,k-1
\end{eqnarray}
(here, $(r,s)=\{x\dvtx f(x)>0\}$; that the support is an interval
follows from (\ref{eq1.2})), where $Q_{i,k}$ are the polynomials defined
by the recurrence above. Obviously, an alternative condition,
sufficient for (\ref{eq3.3}) (and, hence, also for (\ref{eq3.2})), is
%
%
%
\begin{eqnarray}\label{eq3.4}
&&\lim_{x\to r+}q^{i+1}(x)f(x)x^j g^{(i)}(x)=\lim_{x\to
s-}q^{i+1}(x)f(x)x^j g^{(i)}(x)=0,
\nonumber
\\[-8pt]
\\[-8pt]
\nonumber
&&\quad i=0,1,\ldots,k-1, j=0,1,\ldots,k-i-1.
\end{eqnarray}
Assuming that $X$ has $2k$ finite moments, $k\geq1$, it is seen
that (\ref{eq3.4}), and hence (\ref{eq3.2}), is fulfilled for any polynomial $g$
of degree at most $k$. For example, for the upper limit, we have
$\lim_{x\to s-} q(x)f(x)=0$ so that (\ref{eq3.4}) trivially holds if
$s<\infty$; also, if $s=+\infty$, then (\ref{eq3.4}) follows from
$q(x)f(x)=\mathrm{o}(x^{-(2k-1)})$ as $x\to+\infty$, which can be shown
by observing that $q(x)f(x)$ is eventually decreasing, positive
and, by the assumption of finite $2k$th moment, satisfies
$\lim_{x\to+\infty} \int_{x/2}^{x} y^{2k-2}q(y)f(y)\,\mathrm{d}y=0$.
Therefore, the explicit orthogonality relation
is
\begin{equation}\label{eq3.5}
\E[P_k(X)P_m(X)]=\delta_{k,m}
k!\E[q^{k}(X)]\prod_{j=k-1}^{2k-2}(1-j\delta),\qquad
k,m=0,1,\ldots,n,
\end{equation}
where $\delta_{k,m}$ is Kronecker's delta, provided that $X$ has
$2n$ finite moments. The proof follows by a trivial application of
(\ref{eq3.2}) to $g(x)=P_m(x)$, for $0\leq m\leq k\leq n$ (cf. \cite{28}).

It should be noted, however, that the condition (\ref{eq3.4}) or (\ref{eq3.3})
imposes some unnecessary restrictions on $g$. In fact, the
covariance identity (\ref{eq3.2}) (which enables a general form of the
Fourier coefficients of $g$ to be constructed in terms of its
derivatives) holds, in
our case, in its full generality; the proof requires the novel
inversion formula (\ref{eq3.6}) or (\ref{eq3.7}) below, which may be of some
interest in itself.

\begin{theo}\label{th3.1} Assume that $X$ satisfies
\textup{(\ref{eq1.2})} and consider the polynomial $P_k(x)$ defined by \textup{
(\ref{eq1.4})}, where $(r,s)=\{x\dvtx f(x)>0\}$.
\begin{longlist}[(a)]
\item[(a)] If $X$ has $2k-1$ finite moments \textup{($k\geq1$)}, then
the following ``Rodrigues inversion'' formula holds:
%
%
%
\begin{eqnarray}\label{eq3.6}
q^k(x) f(x)
&=&
\frac{(-1)^k}{(k-1)!}\int_{r}^x
(x-y)^{k-1}P_k(y) f(y)\,\mathrm{d}y
\\
\label{eq3.7}
&=&
\frac{1}{(k-1)!}\int_{x}^s
(y-x)^{k-1}P_k(y) f(y)\,\mathrm{d}y,\qquad  x\in(r,s).
\end{eqnarray}
\item[(b)] If $X$ has $2k$ finite moments and
$\E q^{k}(X)|g^{(k)}(X)|<\infty$,
then $\E|P_k(X)g(X)|<\infty$ and
the covariance identity \textup{(\ref{eq3.2})} holds.
\end{longlist}
\end{theo}
\begin{pf}
(a) Let $H_1(x)$, $H_2(x)$ be the left-hand side and right-hand side,
respectively, of
(\ref{eq3.6}). It is easy to see that the integral $H_2(x)$ is finite
(this requires only $2k-1$ finite moments). Moreover, expanding
$(x-y)^{k-1}$ in the integrand of $H_2(x)$ according to Newton's
formula, it follows that $H_2^{(k)}(x)=(-1)^k P_k(x)f(x)$, $x\in
(r,s)$, and, thus, by the definition (\ref{eq1.4}),
$(H_1(x)-H_2(x))^{(k)}=H_1^{(k)}(x)-H_2^{(k)}(x)$ vanishes
identically in $(r,s)$. Therefore, $H_1-H_2$ is a polynomial of
degree at most $k-1$. By (\ref{eq3.1}), $\lim_{x\to r+}H_1^{(i)}(x)=0$ for
all $i=0,1,\ldots,k-1$ because $q(x)f(x)\to0$ as $x\to r+$
and, for the case $r=-\infty$, $q(x)f(x)=\mathrm{o}(x^{-(2k-2)})$ as
$x\to-\infty$ since the $(2k-1)$th moment is finite, $q$ is of
degree at most $2$ and $q(x)f(x)$ is increasing and positive in a
neighborhood of $-\infty$. Similarly, using the fact that $P_k$ is
a polynomial of degree at most $k$ and observing that
$\lim_{x\to r+}x^i \int_{r}^x y^{k+j-1-i} f(y) \,\mathrm{d}y=0$
for all
$i=0,1,\ldots,k-1$, $j=0,1,\ldots,k$ (again, $2k-1$ finite moments
suffice for this conclusion), it follows that
$\lim_{x\to
r+}H_2^{(i)}(x)=0$
for all $i=0,1,\ldots,k-1$. This proves that
$H_1-H_2$ vanishes identically in $(r,s)$ and (\ref{eq3.6}) follows.
Finally, (\ref{eq3.7}) follows from (\ref{eq3.6}) and (\ref{eq3.2}) with
$g(y)=(x-y)^{k-1}$ (the validity of (\ref{eq3.2}) for polynomials $g$
of
degree at most $k-1$ can be shown directly,
using repeated integration by parts, as above).

(b) Suppose that $k\geq1$;
otherwise, since $P_0(x)=1$, we have nothing to show. Since
$\E[P_k(X)]=\E[P_k(X)P_0(X)]=0$ from (\ref{eq3.5}), either $P_k(x)$ vanishes
identically for $x\in(r,s)$ (in which case, (\ref{eq3.2}) trivially holds)
or, otherwise, it must change sign at least once in $(r,s)$.
Assume that $P_k$ is not identically zero and consider the
change-sign points of $P_k$ in $(r,s)$, say,
$\rho_1<\rho_2<\cdots<\rho_m$ (of course, $1\leq m\leq k$
because $P_k$ is a polynomial of degree at most $k$). Fix a point
$\rho$ in the finite interval $[\rho_1,\rho_m]$ and write, with
the help of (\ref{eq3.6}) and (\ref{eq3.7}),
%
%
%
\begin{eqnarray}\label{eq3.8}
&&\E\bigl[q^{k}(X) \bigl|g^{(k)}(X)\bigr|\bigr]\nonumber\\
&&\quad=
\frac{(-1)^k}{(k-1)!}\int_r^{\rho} \bigl|g^{(k)}(x)\bigr|
\int_{r}^x (x-y)^{k-1}P_k(y) f(y)\,\mathrm{d}y \,\mathrm{d}x
\\
&&\qquad{}
+
\frac{1}{(k-1)!}\int_{\rho}^s \bigl|g^{(k)}(x)\bigr| \int_{x}^s
(y-x)^{k-1}P_k(y) f(y)\,\mathrm{d}y \,\mathrm{d}x.\nonumber
\end{eqnarray}
Because of the assumption on $g$, both integrals on the right-hand side of
(\ref{eq3.8}) are finite. We wish to show that we can change the order of
integration in both integrals in the right-hand side of~(\ref{eq3.8}). This will
follow from Fubini's theorem if it can be shown that
%
%
\begin{equation}\label{eq3.9}
I(\rho)=\int_{\rho}^s \bigl|g^{(k)}(x)\bigr| \int_{x}^s
(y-x)^{k-1}|P_k(y)| f(y)\,\mathrm{d}y \,\mathrm{d}x<\infty
\end{equation}
and similarly for the other integral in (\ref{eq3.8}). Since $q(x)f(x)>0$
for all $x\in(r,s)$ (see (\ref{eq1.2})), it follows that $q^k(x)f(x)>0$
for all $x\in(r,s)$ as well. Thus, from (\ref{eq3.7}), we get
%
%
%
\begin{equation}\label{eq3.10}
\int_{x}^s (y-x)^{k-1} P_k(y) f(y)\,\mathrm{d}y >0,\qquad x\in[\rho,s).
\end{equation}
On the other hand, $P_k(x)$ does not change sign in the interval
$(\rho_m, s)$ and, hence, $P_k(x)>0$ for all $x\in(\rho_m, s)$,
showing that
%
%
%
\begin{eqnarray}\label{eq3.11}
&& \int_{\rho_m}^s \bigl|g^{(k)}(x)\bigr| \int_{x}^s
(y-x)^{k-1} |P_k(y)| f(y)\,\mathrm{d}y \,\mathrm{d}x
\nonumber
\\[-8pt]
\\[-8pt]
\nonumber
&&\quad =\int_{\rho_m}^s \bigl|g^{(k)}(x)\bigr| \int_{x}^s
(y-x)^{k-1} P_k(y) f(y)\,\mathrm{d}y \,\mathrm{d}x.
\end{eqnarray}
By the above considerations, it follows that
%
%
%
\begin{eqnarray}
\label{eq3.12} H(\rho)&:=&\inf_{x\in[\rho,\rho_m]}
h(x):=\inf_{x\in[\rho,\rho_m]} \int_{x}^s (y-x)^{k-1} P_k(y)
f(y)\,\mathrm{d}y >0,
\\
\label{eq3.13} S(\rho)&:=&\sup_{x\in[\rho,\rho_m]} s(x)
:=\sup_{x\in[\rho,\rho_m]}\int_x^{\rho_m} (y-x)^{k-1}
|P_k(y)|f(y)\,\mathrm{d}y <\infty,
\\
\label{eq3.14}D(\rho)&:=&\sup_{x\in[\rho,\rho_m]}{d}(x):=\sup_{x\in[\rho,\rho_m]}
\int_{\rho_m}^s (y-x)^{k-1} P_k(y) f(y)\,\mathrm{d}y <\infty
\end{eqnarray}
because the three positive functions $h(x)$, $s(x)$ and ${d}(x)$
defined above are obviously continuous and $x$ lies in the
compact interval $[\rho,\rho_m]$ (note that $h(x)>0$ by (\ref{eq3.10})).
Now, from the inequalities
$s(x)\leq\frac{S(\rho)}{H(\rho)} h(x)$
and ${d}(x)\leq
\frac{D(\rho)}{H(\rho)} h(x)$,
$x\in[ \rho,\rho_m ]$,
we conclude that
%
%
%
\begin{eqnarray}\label{eq3.15}
&&\int_{\rho}^{\rho_m} \bigl|g^{(k)}(x)\bigr| \int_{x}^s
(y-x)^{k-1} |P_k(y)| f(y)\,\mathrm{d}y \,\mathrm{d}x
\nonumber
\\
&&\quad = \int_{\rho}^{\rho_m}
\bigl|g^{(k)}(x)\bigr| \bigl(s(x)+{d}(x)\bigr) \,\mathrm{d}x
\\
&&\quad \leq\frac{S(\rho)+D(\rho)}{H(\rho)}\int_{\rho}^{\rho_m}
\bigl|g^{(k)}(x)\bigr| \int_{x}^s (y-x)^{k-1} P_k(y) f(y)\,\mathrm{d}y
\,\mathrm{d}x.\nonumber
\end{eqnarray}
Combining (\ref{eq3.11}) and (\ref{eq3.15}), we see that there exists a finite
constant $C(\rho)$ (take, for example,
$C(\rho)=\max\{1,(S(\rho)+D(\rho))/H(\rho)\}$) such that
\[
I(\rho)\leq C(\rho)\int_{\rho}^{s}
\bigr|g^{(k)}(x)\bigr| \int_{x}^s (y-x)^{k-1} P_k(y) f(y)\,\mathrm{d}y \,\mathrm{d}x<\infty
\]
and, thus, by (\ref{eq3.9}), we can indeed interchange the order of
integration in the second integral in the right-hand side of (\ref{eq3.8}). Similar
arguments apply to the first integral. By the above arguments and
by interchanging the order of integration in both integrals in the
right-hand side of (\ref{eq3.8})
(with~$g^{(k)}$ in
place of $|g^{(k)}|$), we obtain
%
%
%
\begin{equation}\label{eq3.16}
\E\bigl[q^{k}(X) g^{(k)}(X)\bigr] =\E[P_k(X) G(X)],
\end{equation}
where
%
%
%
\begin{eqnarray}
\label{eq3.17}G(x)
&=&
\frac{(-1)^k}{(k-1)!} \int_x^{\rho} (y-x)^{k-1} g^{(k)}(y)
\,\mathrm{d}y
\\
\label{eq3.18}&=&
\frac{1}{(k-1)!} \int_{\rho}^x (x-y)^{k-1} g^{(k)}(y)
\,\mathrm{d}y,\qquad  x\in(r,s).
\end{eqnarray}
Differentiating (\ref{eq3.17}) or (\ref{eq3.18}) $k$ times, it is easily seen that
$G^{(k)}=g^{(k)}$ so that $G-g=H_{k-1}$ is a polynomial of degree
at most $k-1$ and the desired result follows by (\ref{eq3.16}) and the
orthogonality of $P_k$ to polynomials of degree lower than $k$.
This completes the proof of the theorem.
\end{pf}

\section{An application to lower variance bounds}\label{sec4}

A simple application of Theorem \ref{th2.2} leads to the following lower
variance bound.

\begin{theo}\label{th4.1}
Fix $n\in\{1,2,\ldots\}$ and assume that $X$ satisfies
\textup{(\ref{eq2.1})} and has $2n$ finite moments.
Then, for
any function $g$ satisfying
\begin{equation}\label{eq4.1}
\E  \bigl[
q^{[k]}(X)|\Delta^{k}[g(X)]|
\bigr]<\infty \qquad\mbox{for }
k=0,1,\ldots,n,
\end{equation}
the bound
\begin{equation}\label{eq4.2}
\Var g(X)\geq\sum_{k=1}^n \frac{\E^2[q^{[k]}(X)\Delta
^{k}[g(X)]]}{k!\E
[q^{[k]}(X)]\prod_{j=k-1}^{2k-2}(1-j\delta)}
\end{equation}
holds \textup{(where the $k$th term in the sum should be treated as
zero whenever $\E[q^{[k]}(X)]$ vanishes)} with equality if and only if
$g$ is a polynomial of degree at most $n$.
\end{theo}
\begin{pf} Assume that $\E[ g^2(X)]<\infty$ (otherwise, we
have nothing
to show). By Theorem \ref{th2.1},
the polynomials
$\{P_k/\sqrt{\E[P_k^2(X)]}; k=0,1,\ldots,\min\{n,N\}\}$
form an
orthonormal basis of all polynomials with degree up to $n$, where
$N+1$ is the cardinality of~$J$. Observing that the $k$th Fourier
coefficient for $g$ is, by (\ref{eq2.7}), (\ref{eq2.8}) and Theorem \ref{th2.2},
\[
\frac{\E[
P_k(X)g(X)]}{\E^{1/2}[P_k^2(X)]}
=\frac{\E[q^{[k]}(X)\Delta^{k}[g(X)]]}{(k!\E[q^{[k]}(X)]
\prod_{j=k-1}^{2k-2}(1-j\delta))^{1/2}},\qquad k\leq
\min\{n,N\},
\]
the desired result follows by an application of the finite form of
Bessel's
inequality.
\end{pf}

It is worth mentioning here the similarity of the lower variance
bound (\ref{eq4.2}) with the Poincar\'{e}-type (upper/lower) bound for the
discrete Pearson family, obtained recently in \cite{2},
namely
\begin{equation}\label{eq4.3}
(-1)^{n}\Biggl( \Var g(X)-\sum_{k=1}^{n}
\frac{(-1)^{k+1}}{k!\prod_{j=0}^{k-1}(1-j\delta)}
\E  \bigl[q^{[k]}(X)
(\Delta^{k}[g(X)])^2\bigr]\Biggr)\geq0.
\end{equation}

The following examples can be verified immediately.

\begin{example} \label{ex4.1}If $X$ is $\operatorname{Poisson}(\lambda)$,
then $q(x)=\lambda$ so that $\delta=0$ and (\ref{eq1.6}) follows from
(\ref{eq4.2}) (see also \cite{21}).
Moreover,
the equality in (\ref{eq1.6}) holds if and only if $g$ is a polynomial of
degree at most~$n$.
\end{example}

\begin{example}\label{ex4.2} For the $\operatorname{binomial}(N,p)$
distribution, $q(x)=(N-x)p$ so that $\delta=0$ and $\E
[q^{[k]}(X)]=(N)_k p^k(1-p)^k$. Thus, (\ref{eq4.2}) yields the bound
\[
\Var g(X)\geq\sum_{k=1}^{\min\{n,N\}} \frac{ p^k}{k!
(N)_k(1-p)^k}\E^2\bigl[ (N-X)_k \Delta^{k}[g(X)]\bigr]
\]
with equality only for polynomials of degree at most $n$. Note
that there is equality
if $n\geq
N$.
\end{example}

\begin{example}\label{ex4.3} For the negative
$\operatorname{binomial}(r,p)$ with $p(x)=(r+x-1)_x p^r(1-p)^x/x!$,
$x=0,1,\ldots,$ (\ref{eq2.1}) is satisfied with $q(x)=(1-p)(r+x)/p$. Thus,
$\delta=0$, $\E[q^{[k]}(X)]=(1-p)^k [r]_k /p^{2k}$ and (\ref{eq4.2})
produces the bound
\[
\Var g(X)\geq\sum_{k=1}^n \frac{(1-p)^{k}}{k! [r]_k} \E^2\bigl[
[r+X]_k \Delta^{k}[g(X)]\bigr]
\]
(in the above formulae and elsewhere in the paper,
$[x]_n=x(x+1)\cdots(x+n-1)$ if $n\geq1$ and
$[x]_0=1$).
\end{example}

\begin{example}\label{ex4.4} For the discrete
Uniform$\{1,2,\ldots,N\}$, (\ref{eq2.1}) is satisfied with
$q(x)=x(N-x)/2$; thus, $\delta=-1/2$ and (\ref{eq4.2})
entails the bound
(which is an identity if $n\geq N-1$)
\[
\Var g(X)\geq N \sum_{k=1}^{\min\{n,N-1\}} \frac{(2k+1)(N-k-1)!}{
(k!)^2 (N+k)!}\E^2 \bigl[[X]_k(N-X)_k \Delta^{k}[g(X)]\bigr].
\]
\end{example}

The lower variance bound for the continuous Pearson system is
stated in the following theorem; its proof, being an immediate
consequence of (\ref{eq3.2}), (\ref{eq3.5}), Theorem \ref{th3.1} and a straightforward
application of the finite form of Bessel's inequality
(cf. the proof of Theorem \ref{th4.1} above), is omitted.

\begin{theo}\label{th4.2}
Assume that $X$ satisfies \textup{(\ref{eq1.2})} and has finite moment of
order $2n$ for some fixed $n\geq1$.
Then, for
any function $g$ satisfying $\E[q^k(X)|g^{(k)}(X)|]<\infty$,
$k=0,1,\ldots,n$, we have the inequality
\begin{equation}\label{eq4.4}
\Var g(X)\geq\sum_{k=1}^n \frac{\E^2[q^{k}(X)g^{(k)}(X)]}{k!\E
[q^{k}(X)]\prod_{j=k-1}^{2k-2}(1-j\delta)}
\end{equation}
with equality if and only if $g$ is a polynomial of degree at most
$n$.
\textup{(Note that $\E|X|^{2n}<\infty$ implies
$\delta<(2n-1)^{-1}$ and, thus, that
$\delta\notin\{1,1/2,\ldots,1/(2n-2)\}$ if $n\geq2$.)}
\end{theo}

Some examples now follow.

\begin{example}\label{ex4.5} If $X$ is
$\operatorname{Normal}(\mu,\sigma^2)$, then $q(x)=\sigma^2$ and (\ref{eq4.4}) yields the
Houdr\'{e}--Kagan variance bound (\ref{eq1.7}) (see \cite{20}), under the weaker
assumptions $\E|g^{(k)}(X)|<\infty$, $k=0,1,\ldots,n$; equality
holds if and only if $g$ is a polynomial of degree at most
$n$.
\end{example}

\begin{example}\label{ex4.6} If $X$ is $\Gamma(a,\lambda)$
with density $\lambda^a x^{a-1}\mathrm{e}^{-\lambda x}/\Gamma(a)$, $x>0$,
then $q(x)=x/\lambda$, $\E[q^k(X)]=[a]_k/\lambda^{2k}$ and (4.4)
entails the bound
\[
\Var g(X)\geq\sum_{k=1}^{n} \frac{1}{k![a]_k} \E^2\bigl[ X^k
g^{(k)}(X)\bigr]
\]
with equality only for polynomials of degree at most $n$.
\end{example}

\begin{example}\label{ex4.7} For the standard uniform
density, $q(x)=x(1-x)/2$, $\delta=-1/2$, $\E
[q^{k}(X)]=(k!)^2/(2^k (2k+1)!)$ and we get the bound
\[
\Var g(X)\geq\sum_{k=1}^n \frac{2k+1}{(k!)^2} \E^2\bigl[ X^k (1-X)^k
g^{(k)}(X)\bigr]
\]
with equality only for polynomials of degree at most $n$. Similar
inequalities hold for all beta densities.
\end{example}

Of course, in the above three examples, the
corresponding orthogonal polynomials (\ref{eq1.4}) are the well-known
Hermite, Laguerre and Jacobi (Legendre), respectively, so that one
can alternatively obtain the results using explicit expressions
for the variance and generating functions for the polynomials
(see, e.g., \cite{5,6,12,16,19});
this is also the
case for the discrete polynomials corresponding to Examples
\ref{ex4.1}--\ref{ex4.4} (namely, Charlier, Krawtchouck, Meixner and Hahn,
respectively).
However, as will become clear from the following example,
the considerations given here do not only simplify and unify the
calculations, but also go beyond
the classical
polynomials.

\begin{example}\label{ex4.8} If $X$ follows the $t_N$
distribution (Student's $t$ with $N$ degrees of freedom) with
density
\[
f(x)=\frac{\Gamma((N+1)/2)}{\sqrt{N\uppi}\Gamma(N/2)}
\biggl(1+\frac{x^2}{N}\biggr)^{-(N+1)/2},\qquad x\in\R,
\]
then it is well known that $X$ has only $N-1$ finite integral moments.
However, for $N>1$, $X$~satisfies (\ref{eq1.2}) and its quadratic
$q(x)=(N+x^2)/(N-1)$ has $\delta=1/(N-1)$. Thus, (\ref{eq4.4}) applies for
sufficiently large $N$ (see also \cite{22} for a
Poincar\'{e}-type bound corresponding to (\ref{eq4.3})). To this end, it
suffices to calculate
\[
\E[q^{k}(X)]= \biggl(\frac{N}{N-1}\biggr)^k
\prod_{j=1}^{k}\biggl(1+\frac{1}{N-2j}\biggr),\qquad
 k\leq
(N-1)/2,
\]
and $\prod_{j=k-1}^{2k-2}(1-j\delta)=(N-k)_k/(N-1)^k$.
Theorem \ref{th4.2} yields the (non-classic) bound
\[
\Var g(X)\geq
\sum_{k=1}^{n} \frac{\E^2[ (N+X^2)^k g^{(k)}(X)]}{k!N^k
(N-k)_k\prod_{j=1}^{k}(1+1/(N-2j))},\qquad n\leq\frac{N-1}{2},
\]
with equality only for polynomials of degree at most $n$.
\end{example}

It seems that it would be difficult to work with the explicit forms of
the corresponding orthogonal polynomials, obtained by (\ref{eq1.4}).
Note that a
similar bound can be easily obtained for the Fisher--Snedecor $F_{n_1,n_2}$
distribution and that
Schoutens \cite{31} has obtained the corresponding Stein's equation,
useful in approximating the $t_N$-distribution.

\section{A general variance/covariance representation}\label{sec5}

The main application of the present article (Section \ref{sec4})
presents a convenient procedure for approximating/bounding
the variance of $g(X)$ when $g$ is smooth enough and
when $X$ is ``nice'' enough (Pearson). Although the procedure
is based on the corresponding orthonormal polynomials,
$\phi_k=P_k/\E^{1/2}[P_k^2]$,
the main point is that we do not need explicit forms
for $P_k$. All we need are the Fourier coefficients of $g$,
$c_k=\E[\phi_k g]$, but, due to the identities (\ref{eq3.2}) and (\ref{eq2.7}),
the Fourier coefficients can be simply expressed in terms
of the quadratic $q$ and the derivatives/differences of $g$ when $X$
belongs to the Pearson/Ord family, that is,
when (\ref{eq1.2}) or (\ref{eq2.1}) is satisfied.

A natural question thus arises:
\textit{
is it true that the
$n$th partial sums,
given by the right-hand sides
of \textup{(\ref{eq4.2})} and \textup{(\ref{eq4.4})},
converge to the variance of $g(X)$
as $n\to\infty$}?
More generally, assume that we want to calculate $\Var g(X)$,
when $X$ is an r.v. with finite moments of any order.
Let $\mu$ be
the probability measure of $X$, that is, $\mu$ satisfies
$\mu(-\infty,x]=\Pr(X\leq x)$, $x\in\R$, and assume that the
support of $\mu$, $\operatorname{supp}(\mu)$, is not concentrated on a finite
number of points (otherwise, the following considerations become
trivial). It is well known that there exists an orthonormal
polynomial system (OPS) $\mathcal{F}=\{\phi_0,\phi_1,\ldots\}$,
which can be obtained by an
application of the Gram--Schmidt orthonormalization
process to the
real system $\mathcal{F}_0=\{1,x,x^2,\ldots\}\subset L^2(\R,\mu)$.
Each (real) polynomial $\phi_k$ is of degree (exactly) $k$ and satisfies
the orthonormality condition
$
\E[\phi_i(X)\phi_j(X)]=\int_{\sR} \phi_i(x)\phi_j(x)
\,\mathrm{d}\mu(x)=\delta_{ij}
$.
The members of $\mathcal{F}$ are uniquely determined
by the moments of $X$ (of course, any element of
$\mathcal{F}$ can be multiplied by $\pm1$). For any function
$g\in L^2(\R,\mu)$ (i.e., with finite variance),
we first calculate
the Fourier coefficients
\begin{equation}\label{eq5.1}
c_k=\E[g(X)\phi_k(X)]=\int_{\sR}\phi_k(x)g(x)\,\mathrm{d}\mu(x),\qquad
k=0,1,2,\ldots,
\end{equation}
and then use the well-known Bessel inequality,
\begin{equation}\label{eq5.2}
\Var g(X)\geq\sum_{k=1}^\infty c_k^2,
\end{equation}
to obtain the desirable lower variance bound. Clearly,
Theorems \ref{th4.1} and \ref{th4.2} just provide convenient forms of the
$n$th partial sum in (\ref{eq5.2})
for some particularly interesting cases
(Pearson/Ord system). It is well known
that the sum in the right-hand side
of (\ref{eq5.2}) is equal to the variance (for all $g\in L^2(\R,\mu)$) if and
only if the OPS $\mathcal{F}$ is complete in $L^2(\R,\mu)$ (and, thus, it is
an orthonormal basis of $L^2(\R,\mu)$). This means
that the set of real polynomials,
$\operatorname{span}[\mathcal{F}_0]$,
is dense in $L^2(\R,\mu)$. If this is the case, then
Bessel's inequality, (\ref{eq5.2}), is strengthened to
Parseval's identity,
\begin{equation}\label{eq5.3}
\Var g(X)=\sum_{k=1}^{\infty} c_k^2 \qquad \mbox{for any }g\in
L^2(\R,\mu).
\end{equation}
The following theorem summarizes and unifies
the above observations.
\begin{theo}\label{th5.1}
Assume that $X$ has probability measure $\mu$ and finite moments
of any order.
\begin{longlist}[(a)]
\item[(a)] If $X$ satisfies \textup{(\ref{eq1.2})} and
if $g\in L^2(\R,\mu)\cap D^{\infty}(r,s)$,
then \textup{(\ref{eq5.3})} can
be written as \textup{(cf. Theorem \ref{th4.2})}
\begin{equation}\label{eq5.4}
\Var g(X)= \sum_{k=1}^{\infty} \frac{\E^2[q^{k}(X)g^{(k)}(X)]}
{k!\E[q^{k}(X)]\prod_{j=k-1}^{2k-2}(1-j\delta)},
\end{equation}
provided that the polynomials are dense in
$L^2(\R,\mu)$ and that
\begin{equation}\label{eq5.5}
\E\bigl[q^{k}(X)\bigl|g^{(k)}(X)\bigr|\bigr]<\infty,\qquad k=0,1,2,\ldots.
\end{equation}
\item[(b)] Similarly, if $X$ satisfies \textup{(\ref{eq2.1})}, if $\mu$ is not
concentrated on
a finite integer interval and if $g\in L^2(\R,\mu)$, then
\textup{(\ref{eq5.3})} yields the identity
\textup{(cf. Theorem \ref{th4.1})}
\begin{equation}\label{eq5.6}
\Var g(X)= \sum_{k=1}^{\infty} \frac{\E^2[q^{[k]}(X)\Delta
^k[g(X)]]}{k!\E
[q^{[k]}(X)]\prod_{j=k-1}^{2k-2}(1-j\delta)},
\end{equation}
provided that the polynomials are dense in
$L^2(\R,\mu)$ and that
\begin{equation}\label{eq5.7}
\E \bigl[q^{[k]}(X)|\Delta^k[g(X)]|\bigr]<\infty,\qquad
k=0,1,2,\ldots.
\end{equation}
\end{longlist}
\end{theo}

It should be noted that (\ref{eq5.6}) is always true (and reduces to
a finite sum) if $X$ belongs to the discrete Pearson (Ord) family
and the support of $\mu$ is finite. In this case,
the sum adds zero terms whenever $\E[q^{[k]}(X)]=0$. On
the other hand, it is well known (due to M. Riesz)
that the polynomials are dense in $L^2(\R,\mu)$ whenever $\mu$
is determined by its moments; see \cite{30} or \cite{3}, page 45. An
even simpler sufficient condition is when $\mu$ has a
finite moment generating function at a neighborhood of
zero, that is, when there exists $t_0>0$ such that
\begin{equation}\label{eq5.8}
M_X(t)=\E \mathrm{e}^{tX}<\infty, \qquad  t\in(-t_0,t_0).
\end{equation}
A proof can be found in, for example, \cite{8}.
Since condition (\ref{eq5.8}) can evidently be checked for all Pearson
distributions, we include an alternative proof in the \hyperref[app]{Appendix}.

Taking into account the above, we have the following
covariance representation.

\begin{theo}\label{th5.2}
Assume that $X$ has a finite
moment generating function at a neighborhood of zero.

\begin{longlist}[(a)]
\item[(a)] If $X$ satisfies \textup{(\ref{eq1.2})} and
if $g_i\in L^2(\R,\mu)\cap D^{\infty}(r,s)$, $i=1,2$,
then
\begin{equation}\label{eq5.9}
\Cov[g_1(X),g_2(X)]= \sum_{k=1}^{\infty} \frac{\E[q^{k}(X)g_1^{(k)}(X)]
\cdot\E[q^{k}(X)g_2^{(k)}(X)]}
{k!\E[q^{k}(X)]\prod_{j=k-1}^{2k-2}(1-j\delta)},
\end{equation}
provided that for $i=1,2$,
\begin{equation}\label{eq5.10}
\E\bigl[q^{k}(X)\bigl|g_i^{(k)}(X)\bigr|\bigr]<\infty,\qquad   k=0,1,2,\ldots.
\end{equation}
\item[(b)] Similarly, if $X$ satisfies \textup{(\ref{eq2.1})},
then
\begin{equation}\label{eq5.11}
\Cov[g_1(X),g_2(X)]= \sum_{k=1}^{\infty} \frac{\E[q^{[k]}(X)\Delta
^k[g_1(X)]]
\cdot\E[q^{[k]}(X)\Delta^k[g_2(X)]]}{k!\E
[q^{[k]}(X)]\prod_{j=k-1}^{2k-2}(1-j\delta)},
\end{equation}
where each term with $\E[q^{[k]}(X)]=0$ should be treated as zero,
provided that for $i=1,2$,
\begin{equation}\label{eq5.12}
\E \bigl[q^{[k]}(X)|\Delta^k[g_i(X)]|\bigr]<\infty,\qquad
k=0,1,2,\ldots.
\end{equation}
\end{longlist}
\end{theo}

\begin{pf}
Let $\alpha_k =\E[\phi_k(X)g_1(X)]$ and
$\beta_k =\E[\phi_k(X)g_2(X)]$ be the Fourier coefficients of
$g_1$ and~$g_2$. It is then a standard inner product property in the
Hilbert space that
$\Cov[g_1(X),g_2(X)]=\sum_{k=1}^{\infty} \alpha_k\beta_k$.
Substituting, for example,
\[
\alpha_k=\E[\phi_k(X)g_1(X)]=\frac{\E[P_k(X)g_1(X)]}{\E^{1/2}[P_k^2(X)]}
=\frac{\E[q^{[k]}(X)\Delta^k[g_1(X)]]}{( k!\E[q^{[k]}(X)
\prod_{j=k-1}^{2k-2}(1-j\delta)])^{1/2}}
\]
(and similarly for the continuous case and for $\beta_k$),
we obtain (\ref{eq5.9}) and (\ref{eq5.11}).
\end{pf}

\begin{application}\label{ap5.1} Assume that $X_1,X_2,\ldots,X_{\nu}$
is a
random sample from $\operatorname{Geometric}(\theta)$, \mbox{$0<\theta<1$}, with probability
function
\[
p(x)=\theta(1-\theta)^x,\qquad x=0,1,\ldots,
\]
and let $X=X_1+\cdots+X_{\nu}$ be the complete sufficient
statistic.
The uniformly minimum variance unbiased estimator, UMVUE, of $-\log
(\theta)$
is then given by (see \cite{2})
\[
T_{\nu}=T_{\nu}(X)=\cases{
0, & \quad $\mbox{if }  X=0,$ \vspace*{4pt}\cr
\displaystyle\frac{1}{\nu}+\frac{1}{\nu+1}+\cdots+\frac{1}{\nu+X-1}, & \quad $\mbox{if }
 X\in\{1,2,\ldots\}.$}
\]
Since no simple form exists for the variance of $T_{\nu}$,
the inequalities (\ref{eq4.3}) have been used in \cite{2} in order to prove
asymptotic efficiency.
However, $X$ is negative $\operatorname{binomial}(\nu,\theta)$
and
\[
\Delta^k [T_{\nu}(X)]=\frac{(-1)^{k-1} (k-1)!}{[\nu+X]_k},\qquad   k=1,2,\ldots,
\]
so that one finds from (\ref{eq5.6}) the exact expression (cf. Example \ref{ex4.3})
\[
\Var T_{\nu} = \sum_{k=1}^{\infty} \frac{(1-\theta)^k}{ k^2
\left({\nu+k-1\atop k}\right) }.
\]
Observe that the first term in the series, $(1-\theta)/\nu$, is
the Cram\'{e}r--Rao lower bound.

Now, let $W_{\nu;n}=W_{\nu;n}(X)=[\nu+X]_n/[\nu]_n$ be the UMVUE
of $\theta^{-n}$ ($n=1,2,\ldots$)
and $U_{\nu;n}=U_{\nu;n}(X)=(\nu-1)_n/[\nu-n+X]_n$ be the UMVUE
of $\theta^n$ ($n=1,2,\ldots,\nu-1$). ($W_{\nu;n}(X)$
is a polynomial, of degree $n$, in $X$.)
It follows that
\begin{eqnarray*}
\Delta^k [W_{\nu;n}(X)]
&=&
\cases{
\displaystyle\frac{(n)_k [\nu+X+k]_{n-k}}{[\nu]_n},
&\quad  $k=0,1,\ldots,n,$ \vspace*{2pt}\cr
0, & \quad $k=n+1,n+2,\ldots,$}
\\
\Delta^k [U_{\nu;n}(X)]
&=&
\frac{(-1)^k
[n]_k(\nu-1)_n}{[\nu-n+X]_{n+k}},\qquad
k=0,1,2,\ldots,
\end{eqnarray*}
\vspace{-3pt}
so that
\begin{eqnarray*}
\E\bigl[[\nu+X]_k \Delta^k [W_{\nu;n}(X)]\bigr]
&=&
(n)_k \E
[W_{\nu;n}(X)]=(n)_k \theta^{-n}, \qquad  k=0,1,\ldots,n,
\\
\E\bigl[[\nu+X]_k \Delta^k [U_{\nu;n}(X)]\bigr]
&=&
(-1)^k [n]_k \E
[U_{\nu;n}(X)]=(-1)^k [n]_k \theta^{n}, \qquad  k=0,1,\ldots.
\end{eqnarray*}
Using (\ref{eq2.7}), (\ref{eq2.8}), (\ref{eq5.6}), (\ref{eq5.11}) and Example \ref{ex4.3}, we immediately
obtain the formulae
\begin{eqnarray*}
\Cov[T_{\nu},W_{\nu;n}]
&=&
\theta^{-n}\sum_{k=1}^n
(-1)^{k-1}\frac{(n)_k}{k [\nu]_k}(1-\theta)^k, \qquad n=1,2,\ldots,
\\
\Cov[T_{\nu},U_{\nu;n}]
&=&
-\theta^{n}\sum_{k=1}^{\infty}
\frac{[n]_k}{k [\nu]_k}(1-\theta)^k, \qquad n=1,2,\ldots,\nu-1,
\\
\Cov[W_{\nu;n},W_{\nu;m}]
&=&
\theta^{-n-m}\sum_{k=1}^{\min\{n,m\}}
\frac{(n)_k(m)_k}{k![\nu]_k}(1-\theta)^k,\qquad  n,m=1,2,\ldots,
\\
\Cov[U_{\nu;n},U_{\nu;m}]
&=&
\theta^{n+m}\sum_{k=1}^{\infty}
\frac{[n]_k [m]_k}{k![\nu]_k}(1-\theta)^k, \qquad n,m=1,2,\ldots,\nu-1,
\\
\Cov[W_{\nu;n},U_{\nu;m}]
&=&
\theta^{m-n}\sum_{k=1}^n
(-1)^k\frac{(n)_k [m]_k}{k![\nu]_k}(1-\theta)^k,  \\
&&\hspace*{-68pt}
n=1,2,\ldots,\  m=1,2,\ldots,\nu-1.
\end{eqnarray*}

The above series expansions are in accordance with
the corresponding results on Bhattacharyya bounds
given in \cite{9};
these results are also based
on orthogonality and completeness properties of Bhattacharyya functions,
obtained by Seth \cite{32}.
Similar series expansions
for the variance can be found in \cite{1,25}.
\end{application}

Next, we present a similar application for the exponential
distribution.

\begin{application}\label{ap5.2} Assume that $X_1,X_2,\ldots,X_{\nu}$
is a
random sample from $\operatorname{Exp}(\lambda)$, $\lambda>0$, with density
$f(x)=\lambda \mathrm{e}^{-\lambda x}$, $x>0$,
and let $X=X_1+\cdots+X_{\nu}$ be the complete sufficient
statistic. We wish to obtain the UMVUE of $\log(\lambda)$ and its variance.
Setting $U=\log(X_1)$, we find that
\[
\E U= \int_{0}^{\infty}\mathrm{e}^{-x} \log(x)
\,\mathrm{d}x-\log(\lambda)=-\gamma-\log(\lambda),
\]
where $\gamma=0.5772\ldots$ is Euler's constant.
Therefore, $-\gamma-U$ is unbiased
and it follows
that the UMVUE of $\log(\lambda)$ is of the form
\[
L_{\nu}=L_{\nu}(X)=\E[-\gamma-\log(X_1)|X]=-\log(X)-\gamma+\sum
_{j=1}^{\nu-1}
\frac{1}{j}.
\]
Since $X$ follows a $\Gamma(\nu,\lambda)$ distribution and
$L_{\nu}^{(k)}(X)=(-1)^k (k-1)! X^{-k}$, we obtain from (\ref{eq5.4})
(cf. Example \ref{ex4.6}) the formula
\[
\Var L_{\nu} =\sum_{k=1}^{\infty} \frac{1}{
k^2 \left({\nu+k-1 \atop
k}\right)},
\]
which is quite similar to the formula for $\Var T_{\nu}$ in
Application \ref{ap5.1}. (Once again, the first term in the series, $1/\nu$,
is the Cram\'{e}r--Rao bound.)
Moreover, the series
$\sum_{k=1}^{\infty}
\left(
k^2 \left({\nu+k-1 \atop
k}\right)\right)^{-1}$\vspace{2pt}
can be simplified in a closed form. Indeed, observing that
$\Var L_1=\sum_{k\geq1}1/k^2=\uppi^2/6$ and
taking into account the identity
\[
\frac{1}{k^{2}}\biggl(\pmatrix{\nu+k-1 \cr k}^{-1}-\pmatrix{\nu+k \cr
k}^{-1}\biggr)=\frac{(\nu-1)!}{\nu}
\biggl(\frac{1}{[k]_{\nu}}-\frac{1}{[k+1]_{\nu}}\biggr), \qquad \nu,k=1,2,\ldots,
\]
we have
\[
\Var L_{\nu}-\Var L_{\nu+1}
=\frac{(\nu-1)!}{\nu}
\sum_{k=1}^{\infty}\biggl(\frac{1}{[k]_{\nu}}-\frac{1}{[k+1]_{\nu
}}\biggr)=\frac{1}{\nu^2}
\]
so that
\[
\Var
L_{\nu}=\frac{\uppi^2}{6}-1-\frac{1}{2^2}-\cdots-\frac{1}{(\nu
-1)^2}=\sum
_{k\geq
\nu}\frac{1}{k^2}.
\]
Finally, using the last expression and the fact that
$(k+1)^{-2}<
\int_{k}^{k+1} x^{-2}\,\mathrm{d}x< k^{-2}$, we get the inequalities
$\sum_{k\geq\nu} (k+1)^{-2}=\Var L_{\nu}-\nu^{-2}<\nu^{-1}<\Var
L_{\nu}$, that is,
\[
1< \nu\Var L_{\nu} < 1+1/\nu,
\]
which shows that $L_{\nu}$ is asymptotically efficient.
\end{application}

\begin{appendix}\label{app}
\section*{Appendix}

\textit{A completeness proof under (\textup{\ref{eq5.8})}}.
It is well known that, under (\ref{eq5.8}), $X$ has finite moments of any
order so that the OPS exists (and is unique).
From the general theory
of Hilbert spaces, it is known that $\mathcal{F}$
is a basis (i.e., it is complete) if and only if it is total, that is,
if and only if there does not exist a non-zero function $g\in L^2(\R
,\mu)$
such that $g$ is orthogonal to each $\phi_k$.
Therefore, it suffices to
show that if $g\in L^2(\R,\mu)$ and if
\setcounter{equation}{0}
\begin{equation}\label{eqA.1}
\E[g(X)\phi_k(X)]=0\qquad
\mbox{for all }
k=0,1,\ldots,
\end{equation}
then $\Pr(g(X)=0)=1$.
Since each $\phi_k$ is a polynomial with non-zero leading coefficient,
(\ref{eqA.1})
is equivalent to
\begin{equation}\label{eqA.2}
\E[X^n g(X)]=\int_{\sR}x^n g(x) \,\mathrm{d}\mu(x)=0, \qquad n=0,1,\ldots,
\end{equation}
and it thus suffices to prove that if (\ref{eqA.2}) holds,
then $\E|g(X)|=0$. Since, by assumption, the functions
$x\mapsto x^n$ and $x\mapsto \mathrm{e}^{tx}$ ($|t|<t_0/2$)
belong to $L^2(\R,\mu)$, it follows from the Cauchy--Schwarz inequality
that $\E[\mathrm{e}^{tX}|g(X)|]<\infty$ for $|t|<t_0/2$ and, thus,
that %
\begin{equation}\label{eqA.3}
\E\bigl[\mathrm{e}^{|tX|}|g(X)|\bigr]\leq\E[\mathrm{e}^{-tX}|g(X)|]+\E[\mathrm{e}^{tX}|g(X)|]<\infty,\qquad |t|<t_0/2.
\end{equation}
From Beppo Levi's theorem and (\ref{eqA.3}) it follows that
\[
\sum_{n=0}^{\infty} \E  \biggl[\frac{|tX|^n}{n!}
|g(X)|\biggr]=
\E  \Biggl[ \sum_{n=0}^{\infty} \frac{|tX|^n}{n!}
|g(X)|\Biggr]=
\E\bigl[\mathrm{e}^{|tX|}|g(X)|\bigr]<\infty, \qquad  |t|<t_0/2,
\]
and, therefore, from Fubini's theorem and (\ref{eqA.2}), we get
\begin{equation}\label{eqA.4}
\E[\mathrm{e}^{tX}g(X)]=
\E  \Biggl[ \sum_{n=0}^{\infty} \frac{t^n X^n}{n!}
g(X)\Biggr]=
\sum_{n=0}^{\infty} \frac{t^n}{n!} \E[{X^n}g(X)]=0,
\qquad  |t|<t_0/2.
\end{equation}
Write $g^+(x)=\max\{g(x),0\}$, $g^-(x)=\max\{-g(x),0\}$ and
observe that from (\ref{eqA.2}) with $n=0$,
$\E g(X)=0$. It follows that $\E g^+(X)=\E g^-(X)=\theta$, say,
where $\theta=\E|g(X)|/2$. Clearly, $0\leq\theta<\infty$ because
both $g^+$ and $g^-$ are
dominated by $|g|$ and, by assumption,
$|g|\in L^2(\R,\mu)\subset L^1(\R,\mu)$.
We assume now that
$\theta>0$ and we shall obtain a contradiction.
Under $\theta>0$, we can define two Borel probability measures, $\nu
^+$ and
$\nu^-$, as follows:
\[
\nu^+(A)=\frac{1}{\theta} \int_A g^+(x) \,\mathrm{d}\mu(x),\qquad
\nu^-(A)=\frac{1}{\theta} \int_A g^-(x)\, \mathrm{d}\mu(x), \qquad  A\in\mathcal{B}(\R).
\]
By definition, both $\nu^+$ and $\nu^-$ are absolutely continuous
with respect to $\mu$, with Radon--Nikodym derivatives
\[
\frac{\mathrm{d}\nu^+}{\mathrm{d}\mu}=\frac{1}{\theta} g^+(x),\qquad
\frac{\mathrm{d}\nu^-}{\mathrm{d}\mu}=\frac{1}{\theta} g^-(x),\qquad   x\in\R.
\]
Since, by (\ref{eqA.4}),
\begin{equation}\label{eqA.5}
\int_{\sR} \mathrm{e}^{tx} g^+(x) \,\mathrm{d}\mu(x)= \int_{\sR} \mathrm{e}^{tx}
g^-(x)\,
\mathrm{d}\mu
(x), \qquad |t|<t_0/2,
\end{equation}
it follows that the moment generating functions of $\nu^+$ and $\nu^-$
are finite (for $|t|<t_0/2$) and
identical because, from (\ref{eqA.5}), we have that, for any $t\in(-t_0/2,t_0/2)$,
\[
\int_{\sR} \mathrm{e}^{tx} \,\mathrm{d}\nu^+(x)=
\frac{1}{\theta}\int_{\sR} \mathrm{e}^{tx} g^+(x) \,\mathrm{d}\mu(x)= \frac{1}{\theta}
\int_{\sR} \mathrm{e}^{tx} g^-(x) \,\mathrm{d}\mu(x)=
\int_{\sR} \mathrm{e}^{tx} \,\mathrm{d}\nu^-(x).
\]
Thus, $\nu^+\equiv\nu^-$ and choosing
$A^+=\{x\dvtx g^+(x)>0\}\subseteq\{x\dvtx g^{-}(x)=0\}$,
we are led to the contradiction
$1=\nu^+(A^+)=\nu^-(A^+)=0$.
\end{appendix}

\section*{Acknowledgements} This work was partially supported by the
University of
Athens' Research fund under Grant No.~70/4/5637. We would like to
thank Professor Roger W. Johnson for providing us with a copy of his
paper. Thanks
are also due to an anonymous referee for the careful reading of
the manuscript and for some helpful comments.
N. Papadatos is devoting this work to his six years old
little daughter, Dionyssia, who is suffering from a serious health
disease.

\printhistory

\end{document}